\documentclass[12pt]{article}
\usepackage{amssymb}
\def\eps{\varepsilon} \def\fin{_{\rm fin}} \def\opp{^{\rm opp}}
\def\id{{\rm id}} \def\cir{\,{\scriptstyle \circ}\,}
\def\eins{\hbox{\rm 1}\mskip-4.4mu\hbox{\rm l}}
\def\ap#1{\alpha^1_{#1_1}} \def\am#1{\alpha^2_{#1_2}}
\def\LI#1#2{\Phi^#2_{#1_#2}} \def\RI#1#2{\Psi^#2_{#1_#2}} 
\def\A#1#2{{\alpha^#2_{#1_#2}}} 
\def\CC{{\mathbb C}} \def\TT{{\cal T}}  
 
\def\QED{\hspace*{\fill}$\square$} 
\begin{document}
\parskip3pt \addtolength{\baselineskip}{1pt}

$\quad$\vskip13mm\noindent
{\LARGE \bf New subfactors \\[4pt] 
associated with closed systems of sectors}
\vskip3mm 
\noindent{\large Karl-Henning Rehren}

\noindent{Institut f\"ur Theoretische Physik, 
Universit\"at G\"ottingen (Germany)}

\vskip5mm 
\addtolength{\baselineskip}{-2pt}

\noindent{\small
{\bf Abstract:}
A theorem is derived which (i) provides a new class of subfactors which 
may be interpreted as generalized asymptotic subfactors, and which
(ii) ensures the existence of two-dimensional local quantum field
theories associated with certain modular invariant matrices. }
\vskip5mm\addtolength{\baselineskip}{2pt}

\section{Introduction and results}

We consider type III von Neumann factors throughout. $End\fin(N)$ stands for
the set of unital endomorphisms $\lambda$ with finite dimension
$d(\lambda)$ of a factor $N$.

A {\em closed $N$-system} is a set $\Delta \subset End\fin(N)$ of mutually 
inequivalent 
irreducible endomorphisms such that (i) $\id_N\in\Delta$, (ii) if  
$\lambda\in\Delta$ then there is a conjugate endomorphism 
$\bar\lambda\in\Delta$, and (iii) if $\lambda,\mu\in\Delta$ then 
$\lambda\mu$ belongs to $\Sigma(\Delta)$, the set of endomorphisms which are 
equivalent to finite direct sums of elements from $\Delta$.
 
Let $N\subset M$ be a subfactor of finite index with inclusion homomorphism 
$\iota\in Mor(N,M)$. An {\em extension of the closed $N$-system $\Delta$} 
is a pair $(\iota,\alpha)$, where $\iota$ is as above, and $\alpha$ is
a map $\Delta\to End\fin(M)$, $\lambda \mapsto\alpha_\lambda$, such that 

(E1) $\quad\iota\cir\lambda = \alpha_\lambda\cir\iota$,

(E2) $\quad\iota(Hom(\nu,\lambda\mu))\subset 
Hom(\alpha_\nu,\alpha_\lambda\alpha_\mu)$.

Conditions (E1) and (E2) mean that $(\iota,\alpha)$ is a monoidal
functor from the full monoidal C* subcategory \cite{DR} of
$End\fin(N)$ with objects $\Pi(\Delta)$ (the set of finite products of
elements from $\Delta$) into the monoidal C* category $End\fin(M)$. In
particular, they imply that $\alpha_\lambda$ satisfy the same fusion
rules as $\lambda\in\Delta$, and that $\alpha_{\id_N}=\id_M$ (being an
idempotent within $End\fin(M)$). It follows that if $R_\lambda\in
Hom(\id_N,\bar\lambda\lambda)$ and $\bar R_\lambda\in
Hom(\id_N,\lambda\bar\lambda)$ are a pair of isometries satisfying
the conjugate equations $(1_\lambda\times R_\lambda^*) (\bar R_\lambda\times
1_\lambda)=d(\lambda)^{-1}1_\lambda = (1_{\bar\lambda}\times \bar R_\lambda^*) 
(R_\lambda\times 1_{\bar\lambda})$, and thus implementing left- and
right-inverses $\Phi_\lambda$ and $\Psi_\lambda$ for $\lambda$ (i.e.,
linear mappings which invert the left and right monoidal products with
$1_\lambda$, cf.\ \cite{LRo}), then so do $\iota(R_\lambda)$ and
$\iota(\bar R_\lambda)$ for $\alpha_\lambda$. 
(The notation $\times$ refers to the
monoidal product of intertwiners \cite{DR}.)
In particular
$\alpha_{\bar\lambda}$ is conjugate to $\alpha_\lambda$. 

While $\lambda\in\Delta$ is irreducible by definition, $\alpha_\lambda$ 
may be reducible, and its left- and right-inverses are not unique in
general. But the Lemma below states that the left- and right-inverses
$\Phi_{\alpha_\lambda}$ and $\Psi_{\alpha_\lambda}$ induced by
$\iota(R_\lambda)$ and $\iota(\bar R_\lambda)$ are in fact the unique
standard (minimal) \cite{LRo} ones, provided $\Delta$ is a finite system.  

We state our main result. 

{\bf Theorem: \sl Let $N_1\subset M$ and $N_2\subset M$ be two subfactors 
of $M$, and $(\iota_1,\alpha^1)$ and $(\iota_2,\alpha^2)$ a pair of 
extensions of a finite closed $N_1$-system $\Delta_1$ 
and a finite closed $N_2$-system $\Delta_2$, respectively. Then there 
exists an irreducible subfactor 
$$ A\equiv N_1\otimes N_2\opp \subset B$$
with dual canonical endomorphism
$$ \theta\equiv\bar\iota\cir\iota \simeq
\bigoplus_{\lambda_1\in\Delta_1,\lambda_2\in\Delta_2}
Z_{\lambda_1,\lambda_2}\; \lambda_1\otimes\lambda_2\opp , $$
whose ``coupling matrix'' $Z$ of multiplicities is given by }
$$ Z_{\lambda_1,\lambda_2} = \dim Hom(\ap\lambda,\am\lambda) . $$

Here, $\iota\in Mor(A,B)$ is the inclusion homomorphism with conjugate
$\bar\iota\in Mor(B,A)$.
 
The following special case when $\Delta_i$ are braided systems is of
particular interest for an application in quantum field theory:

{\bf Proposition 1: \sl Assume in addition that the closed systems $\Delta_1$ 
and $\Delta_2$ are braided with unitary braidings $\eps_1$ and
$\eps_2$, respectively, turning $\Pi(\Delta_1)$ and $\Pi(\Delta_2)$
into braided monoidal categories. If for any
$\lambda_i,\mu_i\in\Delta_i$ and any 
$\phi\in Hom(\ap\lambda,\am\lambda)$, $\psi\in Hom(\ap\mu,\am\mu)$,  

\rm (E3) \sl $\quad(\psi\times\phi)\cir \iota_1(\eps_1(\lambda_1,\mu_1)) = 
\iota_2(\eps_2(\lambda_2,\mu_2))\cir (\phi\times\psi$) 

\noindent holds, then the canonical isometry $w_1\in Hom(\theta,\theta^2)$ 
(defined below in the proof of the Theorem) and the braiding operator
$\eps(\theta,\theta)$ naturally induced by the braidings $\eps_1$ 
and $\eps_2\opp$ satisfy}
$$ \eps(\theta,\theta) w_1 = w_1. $$

This result answers an open question in quantum field theory, where possible 
matrices $Z$ are classified which are supposed to describe the 
restriction of a given two-dimensional modular invariant conformal 
quantum field theory to its chiral subtheories, while it is actually
not clear whether any given solution $Z$ does come from a two-dimensional 
quantum field theory. This turns out to be true for a large class of
solutions. 

Namely, let $N_1=N_2=N$ be a local algebra of chiral observables and 
$\Delta_1=\Delta_2=\Delta$ a braided system of DHR endomorphisms.
If the dual canonical endomorphism $\theta_M$ associated with $N\subset M$ 
belongs to $\Sigma(\Delta)$, then $\alpha$-induction \cite{LR,BE}
provides a pair of extensions $(\iota,\alpha^+)$ and $(\iota,\alpha^-)$ 
which satisfies (E1), (E2) as well as (E3) 
\cite[I, Def.\ 3.3, Lemma 3.5 and 3.25]{BE}. The associated coupling 
matrix $Z_{\lambda,\mu}=\dim Hom(\alpha^+_\lambda,\alpha^-_\mu)$ is
automatically a modular invariant \cite{BEK}. By the characterization
of extensions of local quantum field theories given in \cite{LR}, the
subfactor given by the Theorem induces an entire net of subfactors, indexed
by the double-cones of two-dimensional  
Minkowski space. The statement of Proposition 1 is precisely
the criterium given in \cite{LR} for the resulting two-dimensional 
quantum field theory to be local.
Thus, every modular invariant found by the $\alpha$-induction method given in
\cite{BEK} indeed corresponds to a local two-dimensional quantum field theory
extending the given chiral nets of observables. 

The case $N_1=N_2=M$ hence $Z=\eins$ is known for a while \cite{LR},
and was recognized \cite{M} to yield (up to some trivial tensoring with a
type III factor) the type II asymptotic subfactor \cite{O} associated with 
$\sigma(N)\subset N$ where $\sigma\equiv\bigoplus_{\lambda\in\Delta}\lambda$.
As the asymptotic subfactor $M \vee M^c\subset M_\infty$ associated with a 
fixed point inclusion $M^G\subset M$ for an outer action of a group $G$, 
provides the same category of $M_\infty$-$M_\infty$ bimodules as a fixed 
point inclusion for an outer action of the quantum double $D(G)$ on 
$M_\infty$, general asymptotic subfactors in turn are considered 
\cite{O,EK} as generalized quantum doubles. 

Asymptotic subfactors have the properties

(A1) $M\vee M^c \simeq M\otimes M^c$ are in a tensor product position within
$M_\infty$, and every irreducible $M\vee M^c$-$M\vee M^c$ bimodule 
associated with the asymptotic subfactor respects the tensor product, i.e., 
factorizes into an $M$-$M$ bimodule and an $M^c$-$M^c$ bimodule \cite{O}.

(A2) $M$ and $M^c$ are each other's relative commutant in $M_\infty$. We call 
this property of the triple $(M,M^c,M_\infty)$ {\em normality}.

(A3) The system of $M_\infty$-$M_\infty$ bimodules associated with an 
asymptotic subfactor has a non-degenerate braiding \cite{O,I}.

In the type III framework, the analogous property of (A1) is that for 
a subfactor $A\otimes B \subset C$, the dual canonical endomorphism 
$\theta=\bar\iota\cir\iota$ respects the tensor product, i.e., each of
its irreducible components is (equivalent to) a tensor product
$\alpha\otimes\beta$ of endomorphisms of $A$ and $B$, respectively.
We call a subfactor with this property a {\em canonical tensor product 
subfactor (CTPS)} \cite{R1,R2}.

Let $(A,B,C)$ be a joint inclusion of von Neumann algebras, i.e., 
$A\vee B\subset C$. We call $(A,B,C)$ {\em normal}
if $A$ and $B$ are each other's relative commutant in $C$, which is equivalent 
to $A=A^{cc}$ (i.e., $A\subset C$ is normal in standard terminology), and 
$B=A^c$. For $(A,B,C)$ normal, one has $Z(A)=(A\vee B)^c=Z(B)
\supset Z(C)$, so $A$ and likewise $B$ are factors if and only if 
$A\vee B\subset C$ is irreducible, and in this case $C$ necessarily is 
also a factor.

Obviously, the subfactors constructed in the Theorem are CTPS's
(property (A1) of asymptotic subfactors), while we do not know at present 
whether they always share the property (A3) (braiding), which ought to 
be tested with methods as in \cite{I}. Definitely, the joint inclusions 
$(N_1,N_2\opp,B)$ in the Theorem do not share the normality property
(A2) in general. The following Proposition is a characterization of
normality in terms of the coupling matrix, which suggests to
regard normal CTPS's as ``generalized quantum doubles'', beyond
the class of asymptotic subfactors.  

{\bf Proposition 2: \sl Let $A \otimes B \subset C$ be a CTPS of type
III with coupling 
matrix $Z$, i.e., the dual canonical endomorphism is of the form
$$ \theta\simeq \bigoplus_{\alpha\in\Delta_A,\beta\in\Delta_B}
Z_{\alpha,\beta}\;\alpha\otimes\beta , $$
where $\Delta_A\ni\id_A$ and $\Delta_B\ni\id_B$ are two sets of mutually 
inequivalent irreducible endomorphisms in $End\fin(A)$ and
$End\fin(B)$. Then the following conditions are equivalent. 

(N1) The joint inclusion $(A\otimes\eins_B,\eins_A\otimes B,C)$ is
normal, i.e., $A\otimes \eins_B$ and $\eins_A\otimes B$ are each
other's relative commutants in $C$.  

(N2) The coupling matrix couples no non-trivial sector of $A$ to the trivial 
sector of $B$, and vice versa, i.e.,
$$ Z_{\alpha,\id_B}=\delta_{\alpha,\id_A} \quad \hbox{\rm and} \quad
Z_{\id_A,\beta}=\delta_{\beta,\id_B} . $$

(N3) The sets $\Delta_A$ and $\Delta_B$ are closed $A$- and $B$-systems, 
respectively, i.e., they are both closed under conjugation and fusion. 
There is a bijection $\pi:\Delta_A \to\Delta_B$ which preserves the fusion 
rules, i.e.,
$$ \dim Hom(\alpha_1,\alpha_2\alpha_3) = \dim Hom(\pi(\alpha_1),\pi(\alpha_2)
\pi(\alpha_3)) . $$
The matrix $Z$ is the permutation matrix for this bijection, i.e.,}
$$ Z_{\alpha,\beta} = \delta_{\pi(\alpha),\beta} . $$

\section{Indication of Proofs}

For complete proofs, see \cite{R1,R2}.

{\bf Lemma: \sl Let $(\iota,\alpha)$ be an extension of a closed $N$-system 
$\Delta$. Let $R\in Hom(\id_N,\bar\lambda\lambda)$ and 
$\bar R\in Hom(\id_N,\lambda\bar\lambda)$ be a pair of isometries
as before implementing the unique left- and right-inverses \cite{LRo} 
$\Phi_\lambda$ and $\Psi_\lambda$ for $\lambda\in\Delta$. Then 
$\iota(R_\lambda)$ and $\iota(\bar R_\lambda)$ implement left- and 
right-inverses $\Phi_{\alpha_\lambda}$ and $\Psi_{\alpha_\lambda}$ for
$\alpha_\lambda$. 
If $\Delta$ is finite, then $d(\alpha_\lambda)=d(\lambda)$, 
and $\Phi_{\alpha_\lambda}$ and $\Psi_{\alpha_\lambda}$ are the unique
standard left- and right-inverses.}

{\em Proof of the Lemma:} 
The first statement is obvious, since $\iota(R_\lambda)$ and
$\iota(\bar R_\lambda)$ solve the conjugate equations \cite{LRo} for 
$\alpha_\lambda$ if $R_\lambda$ and $\bar R_\lambda$ do so for $\lambda$. If
$\Delta$ is finite, then the minimal dimensions $d(\alpha_\lambda)$
are uniquely determined by the fusion rules, and the latter must coincide
with those of $\lambda\in\Delta$. Hence $d(\alpha_\lambda)=d(\lambda)$. 
Since $d(\lambda)$ are also the dimensions associated with the
pair of isometries $\iota(R_\lambda)$, $\iota(\bar R_\lambda)$, 
the last claim follows by \cite[Thm.\ 3.11]{LRo}. \QED

Thus, general properties of standard left- and right-inverses \cite{LRo} are
applicable. We shall in the sequel repeatedly exploit the trace
property 
$$d(\rho)\Phi_\rho(S^*T)=d(\tau)\Phi_\tau(TS^*)
\qquad {\rm if}\qquad S,T\in Hom(\rho,\tau) $$
for standard left-inverses of $\rho,\tau\in End\fin(M)$,  
their multiplicativity $\Phi_{\rho\tau}=\Phi_\rho\Phi_\tau$, 
as well as the equality of standard left- and right-inverses
$\Psi_\rho=\Phi_\rho$ on $Hom(\rho,\rho)$. 

{\em Proof of the Theorem:}
First notice that the multiplicity of $\id_A$ in $\theta$ is
$Z_{\id_{N_1},\id_{N_2}} = \dim Hom(\id_M,\id_M) = 1$, so the asserted
subfactor is automatically irreducible. 

In order to show that $\theta$ is the dual canonical endomorphism associated 
with a subfactor $A\subset B$, we make use of Longo's characterization 
\cite{L} of canonical endomorphisms in terms of ``canonical triples'' 
(``Q-systems''). It says that $\theta\in End\fin(A)$ is the dual 
canonical endomorphism associated with $A\subset B$ if (and only if)
there is a pair of isometries $w\in Hom(\id_A,\theta)$ and 
$w_1\in Hom(\theta,\theta^2)$ satisfying

(Q1) $\quad w^*w_1 = \theta(w^*)w_1 = d(\theta)^{-1/2} \eins_A$, 

(Q2) $\quad w_1w_1 = \theta(w_1)w_1$, and 

(Q3) $\quad w_1w_1^* = \theta(w_1^*)w_1$. \newpage

In order to construct the Q-system $(\theta,w,w_1)$ in the
present case, we first choose a complete system of mutually
inequivalent isometries  
$W_{(\lambda_1,\lambda_2,l)}\equiv W_l\in A\equiv N\otimes N\opp$,
where $l$ is considered as a multi-index including 
$(\lambda_1\in\Delta_1,\lambda_2\in\Delta_2,l=1,\dots 
Z_{\lambda_1,\lambda_2})$, and put
$$\theta=\sum_l W_l\; (\lambda_1\otimes \lambda_2\opp)(\,\cdot\,)\; W_l^*. $$
The choice of these isometries is immaterial and affects the subfactor
to be constructed only by inner conjugation.

Since $Hom(\id_A,\theta)$ is one-dimensional, the isometry $w$ is already
fixed up to an irrelevant complex phase, and we choose $w=W_0$, where $0$
refers to the multi-index $l=0\equiv(\id_{N_1},\id_{N_2},1)$.
The second isometry, $w_1$, must be of the form
$$ w_1=\sum_{l,m,n} (W_l \times W_m)\cir \TT_{lm}^n\cir W_n^* $$
where $\TT_{lm}^n\in Hom(\nu_1\otimes\nu_2\opp,
(\lambda_1\otimes\lambda_2\opp)\cir(\mu_1\otimes\mu_2\opp))$, since these
operators span $Hom(\theta,\theta^2)$. 

In turn, $\TT_{lm}^n$ must be of the form 
$$ \TT_{lm}^n = \sum_{e_1,e_2} \zeta_{lm,e_1e_2}^n \; T_{e_1} \otimes
(T_{e_2}^*)\opp \qquad (\zeta_{lm,e_1e_2}^n\in\CC) $$
where $T_{e_i}$ constitute orthonormal isometric bases of the intertwiner
spaces $Hom(\nu_i,\lambda_i\mu_i)$, since these operators span
$Hom(\nu_1\otimes\nu_2\opp,
(\lambda_1\otimes\lambda_2\opp)\cir(\mu_1\otimes\mu_2\opp))\equiv
Hom(\nu_1,\lambda_1\mu_1)\otimes Hom(\nu_2\opp,\lambda_2\opp\mu_2\opp)$. 
Note that if $T\in Hom(\alpha,\beta)$ is isometric in $N$, then
$(T^*)\opp \in Hom(\beta,\alpha)\opp \equiv Hom(\alpha\opp,\beta\opp)$ is
isometric in $N\opp$. The labels $e_i$ are again multi-indices of the
form $(\lambda,\mu,\nu,e=1,\dots \dim Hom(\nu,\lambda\mu))$. 

It remains therefore to determine the complex coefficients 
$\zeta_{lm,e_1e_2}^n$, such that $w_1$ is an isometry satisfying Longo's
relations (Q1-3) above. 
To specify the coefficients, we equip the spaces $Hom(\ap\lambda,\am\lambda)$ 
with the non-degenerate scalar products 
$(\phi,\phi'):=\LI\lambda1(\phi^*\phi')$ (where $\LI\lambda i$ stand
for the induced left-inverses for $\A\lambda i$). With
respect to these scalar products, we choose orthonormal bases
$\{\phi_l, l=1,\dots Z_{\lambda_1,\lambda_2}\}$ for all
$\lambda_1,\lambda_2$, and put 
$$ \zeta_{lm,e_1e_2}^n = \sqrt\frac{d(\lambda_2)d(\mu_2)}{d(\theta)d(\nu_2)}
\; \LI\lambda1
[\iota_1(T_{e_1}^*)(\phi_l^* \times \phi_m^*)\iota_2(T_{e_2})\phi_n]. $$

Condition (Q1) is trivially satisfied, since left multiplication of $w_1$
by $w^*$ singles out the term $l=0$ due to $W_0^*W_l=\delta_{l0}$. This
leaves only terms with $\lambda_i=\id_{N_i}$, hence $\mu_i=\nu_i$, for which 
$T_{e_i}$ are trivial and $\sqrt{d(\theta)}\zeta_{0m,e_1e_2}^n=
\delta_{mn}$ (up to cancelling complex phases), so 
$\sqrt{d(\theta)}w^*w_1 = \sum_n W_nW_n^*=\eins_A$. For
$\theta(w^*)w_1$ the argument is essentially the same.

We turn to the conditions (Q2) and (Q3). Whenever we compute either of
the four products occurring, we obtain a Kronecker delta
$W_s^*W_t=\delta_{st}$ for one pair of the labels $l,m,n,\dots$
involved, while the remaining operator parts are of the form 
$$ (W_l \times W_m \times W_k) \left[(T_{e_1} \times 1_{\kappa_1})T_{f_1}
\otimes (((T_{e_2} \times 1_{\kappa_2})T_{f_2})^*)\opp\right] W_n^* , $$
$$ (W_l \times W_m \times W_k) \left[(1_{\lambda_1} \times T_{g_1})T_{h_1}
\otimes (((1_{\lambda_2} \times T_{g_2})T_{h_2})^*)\opp\right] W_n^* $$
for the left- and right-hand side of (Q2), $w_1w_1 = \theta(w_1)w_1$,
and in turn, \newpage
$$ (W_l \times W_m) \left[T_{e_1}T_{f_1}^*
\otimes ((T_{e_2}T_{f_2}^*)^*)\opp\right] (W_n \times W_k)^* , $$
$$ (W_l \times W_m) \left[(1_{\lambda_1}\times
  T_{g_1}^*)(T_{h_1}\times 1_{\kappa_1}) \otimes 
(((1_{\lambda_2}\times T_{g_2}^*)(T_{h_2}\times 1_{\kappa_2}))^*)\opp\right] 
(W_n \times W_k)^*  $$
for the left- and right-hand side of (Q3), $w_1w_1^* = \theta(w_1^*)w_1$.
(In these expressions, we do not specify the respective intertwiner
spaces to which the various operators $T$ belong, since 
these are determined by the context.)

The numerical coefficients multiplying these operators are, respectively,
$$C_{2L}=\sum_s\zeta^s_{lm,e_1e_2}\zeta^n_{sk,f_1f_2},\quad 
C_{2R}=\sum_s\zeta^s_{mk,g_1g_2}\zeta^n_{ls,h_1h_2}$$
for (Q2), and 
$$C_{3L}=\sum_s\zeta^s_{lm,e_1e_2}\overline{\zeta^s_{nk,f_1f_2}},\quad 
C_{3R}=\sum_s\overline{\zeta^m_{sk,g_1g_2}}\zeta^n_{ls,h_1h_2}$$
for (Q3), with a summation over one common label 
$s=1,\dots Z_{\sigma_1,\sigma_2}$ due to the above Kronecker 
$\delta_{st}$ in each case.  

These summations over $s$ can be carried out. Namely,
factors $\zeta^s_{\cdot\cdot,\cdot\cdot}$ are in fact scalar products 
of the form $\LI\sigma1(X\phi_s)=(X^*,\phi_s)$ within
$Hom(\ap\sigma,\am\sigma)$, so summation with the 
operator $\phi_s^*$ contributing to the other factor $\zeta$ yields 
$\sum_s \LI\sigma1(X\phi_s)\phi_s^* = X$. A factor of the form 
$\zeta^\cdot_{\cdot s,\cdot\cdot}$ can also be rewritten with
the help of the trace property for standard left inverses 
as a scalar product $\LI\sigma1(\phi_s^*X)$ within
$Hom(\ap\sigma,\am\sigma)$, and the evaluation of the sum over $s$ is
likewise possible. 

After some transformations, one arrives at
$$C_{2L} \propto
\LI\nu1[\iota_1(T_{f_1}^*(T_{e_1}^*\times 1_{\kappa_1}))
(\phi_l^*\times\phi_m^*\times\phi_k^*)
\iota_2((T_{e_2}\times 1_{\kappa_2})T_{f_2})
\phi_n],$$
$$C_{2R} \propto
\LI\nu1[\iota_1(T_{h_1}^*(1_{\lambda_1}\times T_{g_1}^*))
(\phi_l^*\times\phi_m^*\times\phi_k^*)
\iota_2((1_{\lambda_2}\times T_{g_2})T_{h_2})
\phi_n]
$$
up to a common factor 
$\sqrt\frac{d(\lambda_2)d(\mu_2)d(\kappa_2)}{d(\theta)^2d(\nu_2)}$.
Summing the operators on both sides of (Q2) as above with the
coefficients $C_{2L}$, $C_{2R}$, and noting that the passage from bases 
$(T_{e}\times 1_{\kappa})T_{f}$ to bases $(1_{\lambda}\times T_{g})T_{h}$ 
of $Hom(\nu,\lambda\mu\kappa)$ for any fixed $\nu,\lambda,\mu,\kappa$
is described by unitary matrices, we conclude
equality of both sides of (Q2). 

For (Q3), similar manipulations give
\begin{eqnarray*} C_{3L} \propto \frac{d(\mu_1)d(\mu_2)}{d(\sigma_2)d(\sigma_1)}\;
\LI{\mu_1\lambda}1
[(\phi_l^*\times\phi_m^*)
\iota_2(T_{e_2}T_{f_2}^*)
(\phi_n\times\phi_k)
\iota_1(T_{f_1}T_{e_1}^*)], \qquad\;\;\,\qquad\qquad\\
C_{3R}\propto \frac{d(\sigma_2)d(\sigma_1)}{d(\nu_1)d(\nu_2)}\times
\qquad\qquad\qquad\qquad\qquad\qquad\qquad
\qquad\qquad\qquad\qquad\qquad\qquad \\
\LI{\mu_1\lambda}1
[(\phi_l^*\times\phi_m^*)
\iota_2((1_{\lambda_2}\times T_{g_2}^*)(T_{h_2}\times 1_{\kappa_2}))
(\phi_n\times\phi_k)
\iota_1((T_{h_1}^*\times 1_{\kappa_1})(1_{\lambda_1}\times T_{g_1}))]
\end{eqnarray*}
up to a common factor
$\sqrt\frac{d(\lambda_2)d(\nu_2)d(\kappa_2)}{d(\theta)^2d(\mu_2)}
d(\lambda_1)$.
Summing the operators on both sides of (Q3) as above with the
coefficients $C_{3L}$, $C_{3R}$, and noting that the passage from bases  
$\sqrt\frac{d(\mu)}{d(\sigma)}T_eT_f^*$ to bases 
$\sqrt\frac{d(\sigma)}{d(\nu)}(1_{\lambda}\times
T_{g}^*)(T_h\times 1_{\kappa})$ of $Hom(\nu\kappa,\lambda\mu)$ 
for any fixed $\nu,\kappa,\lambda,\mu$ is again described by a unitary
matrix, we obtain equality of both sides of (Q3).

It remains to show that $w_1$ is an isometry, $w_1^*w_1=1$.

Performing the multiplication $w_1^*w_1$ yields two Kronecker delta's
from the factors $W_l \times W_m$, and two more Kronecker delta's from
the factors $T_{e_1}\otimes (T_{e_2}^*)\opp$. Thus
$$w_1^*w_1=\sum_{ns}\left(\sum_{lm,e_1e_2}
\overline{\zeta^n_{lm,e_1e_2}}\zeta^s_{lm,e_1e_2}\right) W_nW_s^*, $$
and we have to perform the sums over $l,m,e_1,e_2$ (involving, as
sums over multi-indices, the summation over sectors
$\nu_i,\lambda_i,\mu_i\in\Delta_i$, $i=1,2$). 

Again, we rewrite $\zeta^n_{lm,e_1e_2}$ as
a scalar product $(\phi_m,X)$ within $Hom(\ap\mu,\am\mu)$ 
and perform the sum over $m$ similar as before.
In the resulting expression, both sums over $(e_1,\mu_1)$ 
and over $(e_2,\mu_2)$ can be performed after a unitary passage 
from the bases of orthonormal isometries $T_{e}$ of $Hom(\nu,\lambda\mu)$ 
to the bases $\sqrt\frac{d(\lambda)d(\nu)}{d(\mu)}
(1_\lambda\times T_{e'}^*)(\bar R_\lambda\times 1_\nu)$, making use of the 
conjugate equations between $\bar R_\lambda$ (contributing to the new bases)
and $R_\lambda$ (implementing the left-inverses $\Phi_\lambda$ and
hence $\LI\lambda i$). This produces the expression 
$$\sum_{lm,e_1e_2} \overline{\zeta^s_{lm,e_1e_2}}\zeta^n_{lm,e_1e_2} =
\sum_{l,\lambda_1\lambda_2} \frac{d(\lambda_2)^2}{d(\theta)}\; 
\LI\nu1[\RI\lambda2(\phi_l\phi_l^*)\times(\phi_s^*\phi_n)].$$
Here $\RI\lambda2$ ist the standard right-inverse implemented by
$\iota_2(\bar R_{\lambda_2})$ which coincides with $\LI\lambda2$
on $Hom(\A\lambda2,\A\lambda2)$, and can be evaluated by the trace property:
$\RI\lambda2(\phi_l\phi_l^*)= \LI\lambda2(\phi_l\phi_l^*)
=\frac{d(\lambda_1)}{d(\lambda_2)}\LI\lambda1(\phi_l^*\phi_l)=
\frac{d(\lambda_1)}{d(\lambda_2)}$, while the sum over $l$ yields the
multiplicity factor $Z_{\lambda_1,\lambda_2}$. Hence 
$$\sum_{lm,e_1e_2} \overline{\zeta^s_{lm,e_1e_2}}\zeta^n_{lm,e_1e_2} 
  = \left(\sum_{\lambda_1,\lambda_2}
  \frac{d(\lambda_1)d(\lambda_2)Z_{\lambda_1,\lambda_2}}{d(\theta)}\right) 
  \LI\nu1(\phi_s^*\phi_n)=\delta_{sn},$$
and hence $w_1^*w_1=\sum_n W_nW_n^*=1$.

This completes the proof of the Theorem. For the detailed computations, 
cf.\ \cite{R2}. \QED

{\em Proof of Proposition 1:} Left multiplication of $w_1$ with the
induced braiding operator 
$$\eps(\theta,\theta) = \sum_{mlm'l'} (W_{m'}\times W_{l'}) \cir
(\eps_1(\lambda_1,\mu_1) \otimes (\eps_2(\lambda_2,\mu_2)^*)\opp) \cir
(W_l\times W_m)^*
$$
amounts to a unitary passage from bases $T_e\in Hom(\nu,\lambda\mu)$ to
bases $\eps(\lambda,\mu)T_e\in Hom(\nu,\mu\lambda)$. But by (E3), the
coefficients $\zeta^n_{lm,e_1e_2}$ are invariant under these changes
of bases. Hence $\eps(\theta,\theta)w_1=w_1$. \QED
 
{\em Proof of Proposition 2:} The proof is published in 
\cite[Lemma 3.4 and Thm.\ 3.6]{R1}. \QED

\newpage\section{Conclusion}
We have shown the existence of a class of new subfactors associated
with extensions of closed systems of sectors. The proof proceeds by
establishing the corresponding Q-systems in terms of certain matrix
elements for the transition between two extensions. The new subfactors
are canonical tensor product subfactors and include the asymptotic
subfactors. They may be regarded as generalized quantum doubles if
they satisfy a normality condition for which a simple criterium is
given. The new subfactors also include the local subfactors of
two-dimensional conformal quantum field theory associated with certain
modular invariants, thereby establishing the expected existence of
these theories.

\vskip5mm\addtolength{\baselineskip}{-1pt}
\noindent{\large\bf Acknowledgements}
\vskip1mm
I thank Y. Kawahigashi, M. Izumi, T. Matsui, and I. Ojima who made
possible my visit to Japan during summer 1999 where the present work
was completed. I thank notably Y. Kawahigashi for many discussions on
the present construction as well as on ref.\ \cite{BEK}, and H. Kurose
for giving me the opportunity to present these results at the workshop
``Advances in Operator Algebras'' held at RIMS, Kyoto. Financial
support by a Grant-in-Aid for Scientific Research from the Ministry of
Education (Japan) is gratefully acknowledged.

\small\addtolength{\baselineskip}{-.5pt}

\end{document}